%% file: conjugates_descent.tex
\documentclass{article}

\input{preamble}

\begin{document}

\title{On Conjugates and Adjoint Descent}
\author{Asaf Horev and Lior Yanovski}
\maketitle

\begin{abstract}	
  In this note we present an  $\infty$-categorical framework for descent along adjunctions and a general formula for counting conjugates up to equivalence which unifies several known formulae from different fields.
\end{abstract}

\section{Introduction}

The notion of ``conjugate objects'' or ``objects of the same genus'' arises in many fields in mathematics: 
in commutative algebra as objects that become isomorphic after a field extension (\cite{Serre}), in homotopy theory as spaces that have equivalent Postnikov truncations (\cite{Wilkerson}) and in group theory as nilpotent groups that have isomorphic localizations (\cite{DDK}). 
Often, one also has a formula computing the number of conjugates of a given object. 
In the three examples mentioned above, those are in terms of Galois cohomology, $\lim^1$ of a tower of groups and a double coset formula respectively.

The goal of this paper is twofold:
\begin{itemize}
  \item[A.]{
    To unify and generalize the examples above by giving an abstract $\infty$-categorical definition of conjugates (\cref{def:conjugates}) and a general formula for counting them (\cref{thm:main_A}). 
  }
  \item[B.]{
    To prove a descent result which facilitates the construction of the above $\infty$-categorical framework in many cases of interest (\cref{thm:main_B} and its dual, \cref{cor:dual_adj_descent}).
  }
\end{itemize} 

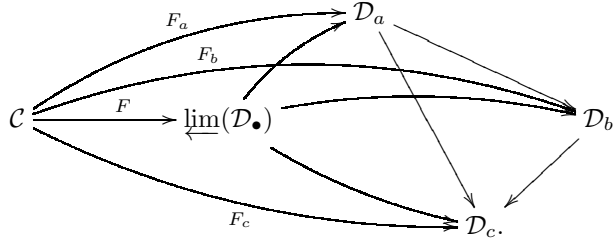
\begin{figure}[!ht]
\begin{align*}
    \vcenter{\xymatrix{
      & & & \D_a\ar[rrd]\ar[rdd] & &\\
      \C\ar@/^1pc/[rrru]^{F_a}\ar@/^1.75pc/[rrrrr]^<<<<<<<<<<<<<{F_b}\ar@/_1pc/[rrrrd]_{F_c}\ar[rr]^F
      & & \invlim(\D_\o)\ar@/^0.5pc/[ru]\ar@/^0.75pc/[rrr]\ar@/_0.5pc/[rrd] & & & \D_b\ar[dl]\\
      & & & & \D_c. &\\
    }}
\end{align*}
\caption{A cone on an $I$-diagram of $\infty$-categories and the comparison functor.}\label{setting_diag}
\end{figure}

For a general definition of conjugate objects, we first need to fix some notation. 
Let $I$ be a simplicial set. 
An $I$-diagram of $\infty$-categories is a map $I\to \Cat$, which we denote by $\D_\o$ (where $\D_a$ is the image of a vertex $a\in I$). 
A cone on $\D_\o$ is an extension of the map $I\to \Cat$ to the cone $I^{\triangleleft}$. 
We denote such a cone by $\C\to\D_\o$, where $\C$ is the image of the cone point. 
In this situation we get a canonical comparison functor $\C\to \invlim (\D_\o)$. 
Now assume that we are in the following setting:
\begin{mysetting} \label{functorial_setting} 
  Let $I$ be a simplicial set, $\D_\o$ an $I$-diagram of $\infty$-categories, and $\C \to \D_\o$ a cone.
  Denote by \(F_a\colon  \C \to \D_a \) the functor corresponding to the edge from the cone point to $a\in I$ and by $F\colon \C \to \invlim (\D_{\o})$ the comparison functor (see \cref{setting_diag}).
\end{mysetting}

Given two objects $x,y\in \C$ one can try to distinguish between them by comparing $F_a(x)$ and $F_a(y)$ in $\D_a$. 
If $F_a(x)$ and $F_a(y)$ fail to be equivalent for some $a\in I$, then clearly $x$ and $y$ cannot be equivalent in $\C$. 
Two objects $x$ and $y$ are called \emph{conjugate} if they can't be distinguished in this way.
More formally, 
\begin{mydef} \label{def:conjugates}
  In the \cref{functorial_setting}, two objects $x$ and $y$ in $\C$ will be called \emph{conjugate} if there exist (not necessarily compatible) equivalences $F_a (x) \simeq F_a(y)$ for every index $a$ in $I$. 
  Let \( \Conj(x) \subseteq \C^\simeq \) denote the full $\infty$-subgroupoid of conjugates of $x$.
\end{mydef}
 
Our main results are:
\begin{theorem}[Conjugates formula] \label{thm:main_A}
  In the \cref{functorial_setting},
  if the comparison functor \( F\colon \C \rightarrow \invlim (\D_\o) \) is an equivalence, then it induces an equivalence of $\infty$-groupoids\footnote{Note that the $\infty$-limit of $\infty$-groupoids corresponds to the homotopy limit of spaces under the identification of the $\infty$-categories of $\infty$-groupoids and spaces.}:
      \begin{align*}
        \Conj(x) \simeq \invlim \BAut F_a(x) 
      \end{align*}
      for every object $x\in \C$. 
      In particular, passing to the set of connected components we get a formula counting the number of conjugates of $x$ up to equivalence:
      \begin{align} \label{eq:conj_formula}
        \pi_0 \Conj(x) \= \pi_0 \invlim \BAut F_a(x). 
      \end{align}
\end{theorem} 
\begin{theorem}[Adjoint descent] \label{thm:main_B}
  In the \cref{functorial_setting},
  if for every $a\in I$ the functor \( F_a\colon \C \to \D_a \) has a right adjoint and $\C$ has all $I$-limits, then the comparison functor 
  \( F\colon \C \to \invlim (\D_\o) \)
  has a right adjoint $G$, given as a composition of two functors:
  \[ \invlim (\D_\o) \xto{G^I} \C^I \xto{\invlim} \C \]
  (\cref{inf_adjoint_descent} below).
  In particular, the comparison functor $F$ is an equivalence if and only if the unit and counit maps are equivalences for all objects.
\end{theorem} 

The proof of \cref{thm:main_A} is completely straightforward (see \cref{sec:conjugates}). 
Its utility comes from the combination with \cref{thm:main_B}. 
To explain this, let us first informally describe the functor $G^I$.
Let $G_a$ be the right adjoint of $F_a$ for each $a\in I$.
An object $y\in \invlim \D_\o$ is roughly a ``compatible system'' of objects \((y_a)_{a\in I}\) with $y_a \in \D_a$.
Using the ``compatibility data'', the collection of objects $G_a(y_a) \in \C$, for each $a\in I$, can be organized in an $I$-diagram, which is defined to be the value of $G^I$ on $y$.
Thus, the functor $G$ itself can be informally described by the following formula:
\begin{align*}
  G((y_a)_{a\in I}) = \invlim_{a\in I} G_a(y_a).
\end{align*}
Moreover, given an object $x\in \C$, the unit of the adjunction is the map
\begin{align*}
	x \to \invlim G_\o (F_\o (x))
\end{align*}
induced by the cone which consists of the units of the individual adjunctions $x \to G_a(F_a(x))$.
Similarly, the $a$-component of the counit, is the map 
\begin{align*}
	F_a(\invlim G_\o (y_\o)) \to y_a
\end{align*}
given by applying $F_a$ to the projection \( \invlim G_\o (y_\o) \to G_a(y_a) \) followed by the counit $F_a(G_a(y_a)) \to y_a$. 

We emphasize that it is \emph{this description} of the adjunction which is the main point of \cref{thm:main_B} and not just the mere existence of the adjoint, which under quite general hypothesis exists by formal arguments (e.g. if all $\infty$-categories are presentable and all functors are left adjoints, this follows from \cite[proposition 5.5.3.13]{Lurie} \footnote{Though there are also some applications for non-presentable $\infty$-categories (such as nilpotent spaces).}). 
One reason for this, is that it provides a concrete formula for the unit and counit, thus reducing the question of whether the comparison functor is fully faithful or an equivalence to an object-wise condition on the unit and counit (as demonstrated in the examples of \cref{sec:examples}). 
Another reason, is that this process can be sometimes applied in reverse.
Given an adjunction, we may get new information about it by representing one of the categories as a limit of a certain diagram and using the specific description of the adjoint provided by \cref{thm:main_B}.
We refer the reader to \cref{Colimits} for a demonstration of this strategy. 

\Cref{thm:main_B} has an obvious dual version obtained by considering the opposites of the $\infty$-categories involved.
Since some of the examples use this dual form, we spell it out for the convenience of the reader.

\begin{cor}[Dual version of \cref{thm:main_B}] \label{cor:dual_adj_descent}
  In the \cref{functorial_setting}, if for every $a\in I$ the functor \( F_a\colon \C \to \D_a \) has a \emph{left} adjoint and $\C$ has all $I^{op}$-colimits, then the comparison functor 
  \( F\colon \C \to \invlim \D_\o \) 
  has a left adjoint $H$ given as a composition
  \begin{align*}
    \invlim (\D_\o) \xto{H^I} \C^{I^{op}} \xto{\colim} \C . 
  \end{align*}
\end{cor}
Denoting by $H_a$ the left adjoint of $F_a$ for all $a\in I$, the functor $H$ can be informally described by the formula
\begin{align*}
	H((y_a)_{a\in I}) = \colim_{I^{op}} H_a(y_a)
\end{align*}
with the analogues description of the unit and counit.

We conclude by briefly sketching the proof of \cref{thm:main_B}.
Since the right adjoint $G$ is constructed as a composition of two functors, $G^I$ and $\invlim$, one might hope to break the problem in two by showing that each of them is a right adjoint. 
While $\invlim$ is of course the right adjoint of the constant diagram functor $\Delta_I$, the functor $G^I$ is usually not the right adjoint of anything. 
To overcome this difficulty we embed $\invlim (\D_\o)$ in a larger category \( \Lax(\D_\o) \), the lax limit of $\D_\o$. 
Viewing the cone $\C \to \D_{\o}$ as a map from the constant diagram on $\C$ to $\D_{\o}$, we show that the induced functor on the \emph{lax limits} 
\begin{align*}
  \C^I\simeq \Lax (\const{\C}) \to \Lax(\D_{\o})
\end{align*}
has a right adjoint, which we denote by $G^I$ (this is an $\infty$-categorical version of a lemma of Barwick (\cite{Barwick2010left}) which is the main ingredient in the proof). 
By composition with the adjunction $\Delta_I\colon \C \adj \C^I\noloc \invlim$, we obtain an adjunction $\C \adj \Lax (\D_{\o})$.
Finally, it remains to observe that the left functor factors through the full subcategory $\invlim \D_{\o}$ and coincides with the comparison functor, hence the adjunction restricts to 
\begin{align*}
	F\colon \C \adj \invlim \D_{\o}\noloc G = \invlim \circ G^I.
\end{align*}

\paragraph{Organization}
In \cref{sec:examples}, we give some examples of conjugates in different settings and sketch the derivation of known counting formulas for them using \cref{thm:main_A} and \cref{thm:main_B}. 
The section concludes with an application of \cref{thm:main_B} to decomposition of colimits.

In \cref{sec:conjugates}, we give a general formula for counting conjugates by proving \cref{thm:main_A}.

In \cref{sec:lax_lim}, we discuss some generalities regarding the lax limit of an $I$-diagram of $\infty$-categories. 
We recall the definition and discuss its functoriality and compatibility with the (non-lax) limit of the diagram.
We also provide some details about the level-wise descriptions of these constructions.

In \cref{sec:lax_adj}, we prove an $\infty$-categorical version of a lemma of Barwick (\cref{lax_adjunction}) showing that a level-wise left adjoint map of $I$-diagrams of $\infty$-categories induces a left adjoint functor on the lax limits with a specific level-wise description.
Finally, we use this to prove \cref{inf_adjoint_descent}, which is a more  detailed version of \cref{thm:main_B}.

\paragraph{Terminology}

We work in the setting of $\infty$-categories \footnote{Also known as 'quasi-categories' or 'weak Kan complexes'.} using heavily the results and terminology of \cite{Lurie}.
In particular, we make extensive use of the technology of Cartesian and coCartesian fibrations.
Therefore, we adopt the notation of \cite{Lurie} whenever possible, unless stated otherwise.
For example, we use the symbols $\invlim$ and $\colim$ for the ($\infty$-categorical) limit and colimit of a functor between $\infty$-categories.
We denote by $\C^\simeq$ the maximal $\infty$-subgroupoid of an $\infty$-category $\C$. 
We abuse notation by identifying ordinary categories with their nerves viewed as $\infty$-categories.

\paragraph{Acknowledgments} We wish to thank our adviser, Emmanuel Farjoun, for suggesting this project to us.
This project is based on his preliminary work \cite{Farjoun}, which inspired the main theorem and examples. 
We also wish to thank Tomer Schlank, Yakov Varshavsky and Assaf Libman for their valuable comments. 

\section{Examples} \label{sec:examples}

In this section we sketch rather informally some examples of the main theorem.
All these examples are well known results and our contribution is only in the unified perspective offered by our main theorems.
The omitted details consist mainly of the formal verification that the respective unit and counit maps of the adjunction produced by our theorem are indeed the expected ones.
The final example is of a somewhat different flavor, as it is about an application solely of \cref{thm:main_B}, highlighting the explicit description of the adjoint.

\begin{ex}[Postnikov towers]
  Following \cite[5.5.6]{Lurie}, let $I$ be the nerve of $\Z_{\geq 0}^{op}$ as a partially ordered set, and $\C$ a presentable $\infty$-category. 
  Let $\tau_{\leq n} \C \subset \C$ denote the full sub-category of $n$-truncated objects.
  The inclusion functor has a left adjoint $\tau_{\le n}\colon \C \to \tau_{\le n}\C$ called the ``$n$-truncation''.
  The truncation functors assemble into an $I$-diagram 
  \[  \cdots \to  \tau_{\le 2}\C \to \tau_{\le 1}\C \to \tau_{\le 0} \C \] 
  that extends to a cone \( \C \to \tau_{\le \o}\C \). 
  Since all the truncation functors \(\tau_{\leq n}\) appearing in the cone are left adjoint to the respective embeddings, we can apply \cref{thm:main_B} to obtain an adjunction \( F\colon \C \adj \invlim \tau_{\leq \o} \C\noloc G \).
  The unit of this adjunction is the canonical map \( X \to \invlim \tau_{\le \o} X \) from an object $X\in \C$ to the limit of its ``Postnikov tower''. 
  The counit for an object $y\in \invlim \tau_{\le \o}\C$ with components $y_n \in \tau_{\le n}\C$ is the natural map $\tau_{\le n}(\invlim y_n) \rightarrow y_n$.
  Hence, if those maps are equivalences for all objects, the functor $F$ is an equivalence and we can apply \cref{thm:main_A}. 
  
  In the classical case where $\C = \mathcal{S}$ is the $\infty$-category of spaces, every space is the homotopy limit of its Postnikov tower, hence the unit of the adjunction is always an equivalence. 
  Similarly, the counit map is always an equivalence since the homotopy groups of the spaces in the tower $(y_n)$ stabilize.  
  This implies the well known fact that the comparison functor 
  $\mathcal{S}\to \invlim \tau_{\le \o}\mathcal{S}$
  is an equivalence.
  Now, by applying \cref{thm:main_A} to a space $X$ we get 
  \[ \pi_0 \Conj(X) \= \pi_0 \invlim \BAut (\tau_{\leq n} X) .\] 
  Finally, the right hand side simplifies by a well known formula (see \cite[lemma 2.2.9]{MayPonto}) to 
  \begin{align*}
    \pi_0 \invlim \BAut (\tau_{\leq n} X) \simeq 
    \lim\nolimits^1 \pi_1 \BAut (\tau_{\leq n} X) \simeq
    \lim\nolimits^1 \hAut (\tau_{\leq n} X) ,
  \end{align*}
  where $\hAut$ is the group of auto-equivalences upto homotopy. 
  Putting everything together, we recover Wilkerson's formula for Postnikov conjugates (see \cite{Wilkerson}):
  \[ \pi_0 \Conj(X) \=  \lim\nolimits^1 \hAut (\tau_{\leq n} X). \]
\end{ex}

\begin{ex}[Galois descent, Borel and Serre \cite{Serre}]
  Let $E/F$ be a finite Galois extension of fields with Galois group $\Gamma$.
  The group $\Gamma$ acts on the category of $E$-vector spaces $\vect_E$ by twisting the coefficients.
  Identifying $\Gamma$ with a single object category, this action turns $\vect_E$ into a coherent $\Gamma$-diagram.
  The extension of scalars functor \( \vect_F \to \vect_E \) transforms coherently under the action of $\Gamma$, defining a coherent $\Gamma$-cone, with homotopy limit \( (\vect_E)^{h\Gamma} \) called the homotopy fixed points.
  Since the extension of scalars functor is a left adjoint, \cref{thm:main_B} produces an adjunction \( \vect_F \leftrightarrows (\vect_E)^{h\Gamma} \) with right adjoint given on an $E$-vector space $W$ by taking the $\Gamma$-fixed points of the underlying $F$-vector space $W_F$.
  The unit and counit of this adjunction are the obvious maps
  \begin{align*}
    \eta\colon V \to\left(E\otimes V\right)_{F}^{\Gamma},\quad 
    \varepsilon\colon E \otimes\left(W_{F}\right)^{\Gamma}\to W .
  \end{align*}
  The Galois property implies $\eta$ is a natural isomorphism, and classical ``Galois descent'' implies $\varepsilon$ is a natural isomorphism as well.
  Hence, we obtain an equivalence of categories \( \vect_F \simeq (\vect_E)^{h\Gamma} \).

  Similarly, and more interestingly, we can consider the categories $\C_F$ and $\C_E$ of vector spaces equipped with tensors of specified signature satisfying some identities (e.g. associative algebras, Lie algebras, Frobenious algebras, quadratic spaces etc.).
  A similar argument gives an equivalence \( \C_F \simeq (\C_E)^{h\Gamma} \)
  \footnote{This requires showing equivalence between suitably constructed larger categories and restricting to a proper subcategory as in \cref{PrinBundles}. We leave the details to the reader.}.
  In all such cases the conjugates formula \eqref{eq:conj_formula} of \cref{thm:main_A} specializes to non-abelian group cohomology
  \[
    \pi_0 \Conj(X) = \pi_0 \BAut(E \otimes X)^{h\Gamma} = H^1 (\Gamma ; \Aut(E \otimes X)) .
  \]
  Thus, recovering the classical formula. 
\end{ex}

\begin{ex}[Arithmetic square for nilpotent spaces, \cite{DDK}]
  In this example we rely on the excellent book \cite{MayPonto} for the relevant facts about localizations of nilpotent spaces. 
  In our framework, we can summarize the situation as follows. 
  Let $\mathbf{Nil}$ be the full subcategory of the $\infty$-category of spaces spanned by nilpotent spaces.
  More generally, given a set of prime numbers  $T$, let \(\mathbf{Nil}_T\) be the full subcategory spanned by $T$-local nilpotent spaces.
  The existence of $T$-localizations and their universal property implies that the inclusion $\mathbf{Nil}_T \hookrightarrow \mathbf{Nil}$ has a left adjoint, taking each space $X$ to its $T$-localization $X_T$. 
  
  Let $T_i$ be a set of prime numbers for each $i=1,\ldots,n$. 
  Denote their union by \( T = \bigcup T_i\) and intersection by \( S = \bigcap T_i\). 
  Assume further that $T_i \cap T_j =S$ for all $i \ne j$ (e.g. $T_1=\{p\}$ and $T_2$ the set of all primes which are not $p$. 
  In this case, $T$ is the set of all primes and $S$ is empty. 
  Hence, $T$-localization is the identity functor and $S$-localization is rationalization).

  The collection of localization functors $\mathbf{Nil}_{T_i}\to \mathbf{Nil}_S$ forms a multi-span diagram of $\infty$-categories
  \begin{align*}
    \vcenter
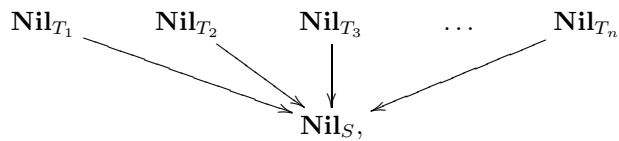
{\xymatrix{
      \mathbf{Nil}_{T_1} \ar[drr] & \mathbf{Nil}_{T_2}\ar[dr] & \mathbf{Nil}_{T_3}\ar[d] & \ldots & \mathbf{Nil}_{T_n}\ar[dll]\\
      & & \mathbf{Nil}_{S},
    }}
  \end{align*}
  which we denote by $\mathbf{Nil}_{\o}$ and the collection of localization functors $\mathbf{Nil}_T\to \mathbf{Nil}_{T_i}$ together with $\mathbf{Nil}_T\to \mathbf{Nil}_{S}$ forms a cone on this diagram. 

  Since all these functors are left adjoints by construction, by \cref{thm:main_B} we obtain an adjunction $\mathbf{Nil}_T \adj \invlim \mathbf{Nil}_{\o}$. 
  The component of the unit of this adjunction at a space $X$ is the map from $X$ to the limit of the diagram of spaces:
  \begin{align*}
    \vcenter
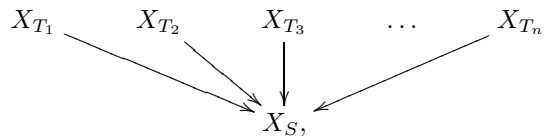
{\xymatrix{
      X_{T_1} \ar[drr] & X_{T_2}\ar[dr] & X_{T_3}\ar[d] & \ldots & X_{T_n}\ar[dll]\\
      & & X_{S} ,
    }}
  \end{align*}
  where the maps are the localization maps. 
  Using the standard procedure of expressing a multi-pullback as an ordinary pullback, this is the same as the map from the top left corner to the pullback int the following square:
  \begin{align*}
    \vcenter{\xymatrix{
      X \ar[d] \ar[r] & \prod\limits_{i=1}^n X_{T_i} \ar[d] \\
      X_S \ar[r]^-{\Delta} & \prod\limits_{i=1}^n X_S .
    }}
  \end{align*}
  
  Since this square is known to be a pullback square (see \cite[8.1.3 and 8.1.4]{MayPonto}), the unit map is an equivalence. 
  Similarly, an object of $\invlim \mathbf{Nil}_{\o}$ consists of spaces $X_i\in \mathbf{Nil}_{T_i}$, $Y\in \mathbf{Nil}_S$ and equivalences  $(X_i)_S \overset{\sim}{\to} Y$. 
  Its image in $\mathbf{Nil}_T$ is a space $Q$ fitting into a pullback diagram 
  \begin{align*}
    \vcenter{\xymatrix{
      Q \ar[d] \ar[r] & \prod\limits_{i=1}^n X_{i} \ar[d] \\
      Y \ar[r]^-{\Delta} & \prod\limits_{i=1}^n Y .
    }}
  \end{align*}
  To verify that the counit map is an equivalence, it is enough to show that the induced localization maps $Q_{T_i}\to X_i$ are equivalences for all $i=1,\ldots,n$, which is indeed true (see \cite[8.1.7 and 8.1.10]{MayPonto}).
    
  We can now apply \cref{thm:main_A} to get a formula for the conjugates of a $T$-local nilpotent space $X$. 
  Unwinding the definitions we get:
  \begin{align*}
    \pi_0 \Conj(X) =  \hAut(X_S) \backslash 
    \prod_{i=1}^n \hAut(X_S) / 
    \prod_{i=1}^n \hAut(X_{T_i}).
  \end{align*}
  This formula is a variant of the known double coset formula for the Mislin genus (see \cite[p. 168]{MayPonto}).
  A similar (and simpler) analysis applies to the special case of nilpotent groups. 
  We refer the reader to \cite[section 7]{MayPonto} for details.
\end{ex}

\begin{ex}[Principal bundles and \v{C}ech cohomology] \label{PrinBundles}
  In this example we work with the 1-category of topological spaces $\mathbf{Top}$.
  Fix a compact topological group $G$ (e.g. $O(n)$) and for each space $X$ consider $\mathbf{Bun}_G(X)$, the groupoid of principal $G$-bundles over $X$.
  Given an open cover $X= \bigcup U_\alpha$ of $X$, we construct a map \( \U = \coprod U_\alpha \to X \) and the associated \v{C}ech nerve \( \U_\o \colon  \Delta^{op} \to \mathbf{Top}\) with \( \U_n = \U \times_X \cdots \times_X \U\). 
  Since principal $G$-bundles pull back along continuous maps, we get an augmented cosimplicial diagram of groupoids \(\mathbf{Bun}_G(X) \to \mathbf{Bun}_G(U_\o)\), which induces a map $\mathbf{Bun}_G(X) \to \invlim \mathbf{Bun}_G(U_\o)$. 
  We would like to demonstrate the well known fact that this map is an equivalence, but considering this as a diagram of categories does not seem to fit into our framework, since maps of groupoids have no adjoints unless they are equivalences (and they aren't in this case). 
  Nevertheless, it is possible to enlarge these categories and functors in a way which will fit the hypothesis of \cref{thm:main_B}. 
  
  Let $\mathbf{Top}_G(X)$ be the category of $G$-spaces with an equivariant map to $X$ considered as a $G$-space with the trivial $G$-action. 
  We will denote such an object by $Y\to X$ and refer to it as a $G$-space over $X$.
  We have a fully faithful embedding $\mathbf{Bun}_G(X) \hookrightarrow \mathbf{Top}_G(X)$. 
  Moreover, a continuous map $f\colon X\to Y$ induces a pullback functor $f^*\colon \mathbf{Top}_G(Y)\to \mathbf{Top}_G(X)$ extending the pullback functor of principal bundles. 
  Only this time, it has a  \emph{left} adjoint $f_!\colon  \mathbf{Top}_G(X)\to \mathbf{Top}_G(Y)$ given simply by post composition with $f$. 
  Applying \cref{cor:dual_adj_descent} we obtain an adjunction:
  \begin{align*}
    R\colon \mathbf{Top}_G(X) \adj \invlim \mathbf{Top}_G(U_\o)\colon L.
  \end{align*}
  Given $Y\to X$ in $\mathbf{Top}_G(X)$, let us denote the counit of the adjunction induced by the map $\U_n \to X$ by $Y|_{\U_n} \to Y$. 
  The counit of the adjunction $R \vdash L$ at the object $Y \to X$ is the map $\colim Y|_{U_{\o}} \to Y$ and is well known to be an isomorphism (see \cite[proposition 4.1]{Vistoli}).
  Similar considerations show that the unit of the adjunction is an isomorphism as well. 
  Hence, we obtain an equivalence of categories $\mathbf{Top}_G(X) \iso \invlim \mathbf{Top}_G(U_\o)$. 
  Since being a principal $G$-bundle is a local condition, this restricts to an equivalence of groupoids $\mathbf{Bun}_G(X) \iso \invlim \mathbf{Bun}_G(U_\o)$.
  
  Applying the formula for conjugates to the trivial $G$-bundle $\tau_X$ on $X$ we get 
  \begin{align*}
    \pi_0(\Conj(\tau_X)) \simeq 
    \pi_0(\invlim_{\Delta} \BAut(\tau_{U_\o})) \simeq
    \pi_0(\invlim_{\Delta} \BMap(U_{\o},G)).
  \end{align*}
  The left hand side is the set of isomorphism classes of principal $G$-bundles trivialized by the cover $\U$. 
  For the right hand side, note that since we are dealing with 1-truncated spaces, we can truncate the diagram at level $2$ without changing its limit. 
  Namely, we can equivalently consider the limit of the truncated cosimplicial diagram (we omitted the degeneracies for typographical reasons):
  \begin{align*}
    \xymatrix{
    \BMap(\U,G) \ar@<.5ex>[r] \ar@<-.5ex>[r] &
    \BMap(\U\times_X \U,G) \ar[r] \ar@<1ex>[r] \ar@<-1ex>[r] &
    \BMap(\U\times_X \U \times_X \U,G)}.
  \end{align*}
  Unwinding the definitions, we see that $\pi_0$ of this limit can be identified with collections of continuous functions $g^{\alpha}_{\beta}\colon U_{\alpha} \cap U_{\beta} \to G$ for all pairs of indexes $\alpha$ and $\beta$ which satisfy the cocycle condition 
  $g^{\beta}_{\gamma} \cdot g^{\alpha}_{\beta} = g^{\alpha}_{\gamma}$ for triples of indexes $\alpha$, $\beta$ and $\gamma$ over the triple intersection $U_{\alpha} \cap U_{\beta} \cap U_{\gamma}$, modulo the equivalence relation $(g^{\alpha}_{\beta})\sim (h^{\alpha}_{\beta})$ if there exists a collection of maps $f_{\alpha}\colon U_{\alpha}\to G$ such that 
  $g^{\alpha}_{\beta}\cdot f_{\alpha} = f_{\beta}\cdot h^{\alpha}_{\beta}$. 
  In other words, we recover the usual formula that identifies the isomorphism types of principal $G$-bundles trivialized by $\U$ with the the \v{C}ech cohomology set $\check{H}^1_{\U}(X; G)$.
\end{ex} 

\begin{ex}[Colimits and decomposition of diagrams] \label{Colimits}
  This example is closely related, and should be compared with, the material discussed in \cite[section 4.2.3]{Lurie}.
  Let $f\colon K\to \C$ be a $K$-diagram in an $\infty$-category $\C$.
  It is possible to compute the colimit of $f$ by breaking $K$ into smaller pieces, computing the colimit over each piece separately, and then taking the colimit of all these partial colimits. 

  More formally, let $K_{\o}\colon I\to \Cat$ be an $I$-diagram of $\infty$-categories with colimit $K$. 
  This induces an $I^{op}$-diagram of functor categories $\C^{K_\o}$ with limit $\C^K$.
  The functors $F_a\colon \C^K \to \C^{K_a}$ take each $K$-diagram $f$ to its restriction along $K_a \to K$, which we denote by $f|_{K_a}$. 
  Consider the constant diagram functor $F\colon \C \to \C^{K}$.
  $F$ has a left adjoint if and only if $\C$ has all $K$-indexed colimits, in which case the left adjoint is the functor $\colim_K$. 
  On the other hand, we can represent $F$ as the comparison functor induced by the cone $\C \to \C^{K_{\o}}$ with all $F_a$ being the constant diagram functors. 
  If $\C$ has all $K_a$-indexed colimits for all $a\in I$, each of the functors $F_a$ admits a left adjoint given by $\colim_{K_a}$.
  If in addition $\C$ has all $I$-colimits, we can apply \cref{cor:dual_adj_descent} to construct a left adjoint to $F$.
  Namely, we conclude that each $K$-diagram $f\colon K\to \C$ has a colimit given by the formula:
  \begin{align*}
    \colim_K f = \colim_{a\in I} (\colim_{K_a} f|_{K_a}).
  \end{align*}
  This in turn can serve as an alternative argument for some of the results in \cite[section 4.4.2]{Lurie}, such as propositions 4.4.2.2, 4.4.2.4 and 4.4.2.6.
\end{ex}

\section{Counting Conjugates} \label{sec:conjugates}
In this section we prove \cref{thm:main_A}.
The equivalence $\C \simeq \invlim (\D_\o)$ implies an equivalence of maximal $\infty$-subgroupoids (i.e. Kan complexes) \( \C^\simeq \simeq (\invlim \D_\o)^\simeq \).
Since the maximal $\infty$-subgroupoid functor is a right adjoint by \cite[proposition 5.2.3.5]{Lurie}, we can identify \((\invlim \D_\o)^\simeq \) with \(\invlim (\D_\o^\simeq) \).

\begin{notation}
  Let $x$ be an object of $\C$.
  Denote by \( \D_a^\simeq(x) \subset \D_a^\simeq \) the connected component of the full $\infty$-subgroupoid of $\D_a^\simeq$ containing \(F_a(x) \in \D_a^\simeq \). 
\end{notation}

Observe that each $\D_a^\simeq(x)$ is a connected $\infty$-groupoid with a distinguished point $F_a(x)$. 
Hence we obtain a canonical equivalence 
$\D_a^\simeq \simeq \BAut F_a(x)$. 
We are therefore reduced to proving the following lemma.

\begin{lem}
  There is a canonical equivalence \( \Conj(x) \simeq \invlim \D_\o^\simeq(x) \). 
\end{lem}

\begin{proof}
  Applying the functor \( \pi_0 \) to the $I$-diagram \(D_\o^\simeq \) we obtain an $I$-diagram of discrete spaces \( \pi_0 (\D_\o)^\simeq \) and a natural transformation of $I$-diagrams
  \( \D_\o^\simeq \to \pi_0 (\D_\o)^\simeq \).
  On the other hand, let $\Delta^0_\o$ be the constant $I$-diagram on  $\Delta^0$. 
  There is a map of $I$-diagrams 
  $\Delta^0_\o\to \pi_0 (\D_\o^\simeq) $
  that takes $\Delta^0_a$ to the connected component of $F_a(x)$. 
  By definition, we have a pullback square
  \begin{align*}
    \vcenter{\xymatrix{
      \Conj(x) \pullbackcorner \ar[d] \ar[rr] & & **[r] \invlim \D_\o^\simeq \simeq \C^\simeq \ar[d] \\
      \Delta^{0} = \invlim \Delta^0_\o \ar[r] & \invlim \D_\o^\simeq \ar[r] &  \invlim \pi_0 (\D_\o^\simeq ).
    }}
  \end{align*}  
  On the other hand, we observe that for each $a\in I$, we have a pullback diagram 
   \begin{align*}
    \vcenter{\xymatrix{
      \D_a^\simeq(x) \pullbackcorner \ar[d] \ar[rr] & & \D_a^\simeq \ar[d] \\
      \Delta^0 \ar[r]^{F_a(x)} & \D_a^\simeq \ar[r] & \pi_0 (\D_a^\simeq).
    }}
  \end{align*}
  Since $I$-limits commute with pullbacks, the lemma follows.
\end{proof}

\section{Lax Limits of $\infty$-Categories} \label{sec:lax_lim}
In this section we discuss (in a rather ad hoc fashion) the lax limit of a diagram of $\infty$-categories, its functoriality and its relation to the ordinary ($\infty$-categorical) limit. 
All proofs are by more or less straightforward manipulations of coCartesian fibrations.

We start by recalling the definition of a lax limit of an $I$-diagram of $\infty$-categories. 
\begin{mydef} 
  Let $I$ be a simplicial set and \( \C_\o\colon  I \to \Cat \) an $I$-diagram of $\infty$-categories, and $q\colon \tint \C_\o \to I$ the coCartesian fibration classifying $\C_\o$.
  We define the \emph{lax limit} of $\C_\o$ to be the $\infty$-category of sections of this coCartesian fibration and denote it by \(  \Lax(\C_\o) = \Map_I (I, \tint \C_\o ) \).
  Namely, $\Lax(\C_\o)$ is the pullback of the following diagram:
  \begin{align} \label{diag:Lax_as_PB}
    \vcenter{\xymatrix{
      \Lax(\C_\o) \ar[r] \ar[d] \pullbackcorner & (\tint \C_\o)^I \ar[d]^{q^I} \\
      \Delta^0 \ar[r] & I^I    
     }}
  \end{align}
where the lower horizontal map picks the identity of $I$. According to \cite[proposition 3.3.3.2]{Lurie} we can identify $\invlim(\C_\o)$ with the full subcategory of $\Lax(\C_\o)$ spanned by the coCartesian sections.
\end{mydef}

\begin{ex}\label{Const_diag}
  For an $\infty$-category $\C$ denote by $\const{\C}$ the constant $I$-diagram on $\C$.
  The coCartesian fibration classified by $\const{\C}$ is the projection $\C\times I \to I$.
  It follows that $\Lax (\const{\C}) = \C^I$.
\end{ex}

\begin{rem}
  We can think of an object \( x\in \Lax(\C_\o) \) given by a section \(x\colon I \to \tint \C_\o \) as a ``lax compatible system'' of objects \( (x_a)_{a\in I} \), where $x_a \in \C_a$ is the value of the section $x$ on the vertex $a\in I$.
  The values of the section $x$ on higher simlicies of $I$ encode the (lax) compatibility of the system.
\end{rem}

\begin{rem}
  The notation \( \Lax(\C_\o) \) is not used in \cite{Lurie}.
  An alternative construction of the lax limit is given in \cite{GHN} and shown to be equivalent to the $\infty$-category of sections.
\end{rem}

We now discuss the functoriality of the lax limit.

\begin{mydef}
  A map of $I$-diagrams of $\infty$-categories from $\C_\o$ to $\D_\o$ is a map \( I \times \Delta^1 \to \Cat \), which restricts to $\C_\o$ and $\D_\o$ on $I\times \Delta^{\set{0}}$ and $I \times \Delta^{\set{1}}$ respectively.
  We denote such a map by \( F_\o\colon \C_\o \to \D_\o \) and for  each  $a\in I$, by $F_a\colon \C_a \to \D_a$ the functor corresponding to the restriction of $F_\o$ to $\set{a} \times \Delta^1$.
\end{mydef}

A map of $I$-diagrams of $\infty$-categories induces a functor between their lax limits in the following way.

\begin{construction} \label{const:assoc_lax}
A map of $I$-diagrams of $\infty$-categories \( F_\o\colon \C_\o \to \D_\o \) classifies a coCartesian fibration \( p\colon  \tint F_\o \to I\times \Delta^1 \).
  Let \( \iota\colon  \Delta^1 \to (I \times \Delta^1 )^I \) be the map corresponding to the identity \( I \times \Delta^1 \to I \times \Delta^1 \). 
  Denote the pullback of \( p^I\colon  (\tint F_\o)^I \to (I\times \Delta^1)^I \) along $\iota$ by \( \Lax(p)\colon  \E \to \Delta^1 \):
  \begin{align} \label{diag:assoc_lax_PB}
    \vcenter{\xymatrix{
      \E \ar[d]_{\Lax(p)} \ar[r] \pullbackcorner & (\tint F_\o)^I \ar[d]^{p^I} \\
      \Delta^1 \ar[r]^-{\iota} & (I \times \Delta^1 )^I.
    }}
  \end{align}
\end{construction}

\begin{lem} \label{lem:assoc_lax}
  Let \(F_\o\colon  \C_\o \to \D_\o \) be a map of $I$-diagrams of $\infty$-categories classifying a coCartesian fibration \( p\colon \tint F_\o \to I\times \Delta^1 \). 
  The morphism \( \Lax(p)\colon  \E \to \Delta^1 \) of \cref{const:assoc_lax} is a coCartesian fibration classified by a functor \( \Lax(F_\o)\colon \Lax(\C_\o) \to \Lax(\D_\o) \).
\end{lem}

\begin{proof}
  By the dual of \cite[proposition 3.1.2.1]{Lurie} the map \( p^I \colon  (\tint F_\o)^I \to (I\times \Delta^1)^I \) is a coCartesian fibration.
  By \cite[proposition 2.4.2.3]{Lurie} the map $\Lax(p)$ is a coCartesian fibration, hence classifies a functor \( \Lax(F) \colon  \E_0 \to \E_1 \) between the fibers of \( \Lax(p) \). 

  We verify that the fiber $\E_0$ of \(\Lax(p)\colon \E \to \Delta^1 \) is isomorphic to $\Lax(\C_\o)$ (proving that $\E_1$ is isomorphic to \(\Lax(\D_\o)\) is similar).
  Factoring $\C_\o$ via $F_\o$ and pulling along the factorization we have a pullback diagram (see \cite[Section 3.3.2]{Lurie})
  \begin{align*}
    \vcenter{\xymatrix{
      \tint \C_\o \ar[d]_{p_0} \ar[r] \pullbackcorner & \tint F_\o \ar[d]_{p} \\
      I\times \Delta^{\set{0}}  \ar[r] \ar@/_1pc/[rr]_{\C_\o} & I \times \Delta^1 \ar[r]^{F_\o} & \Cat,
    }}
  \end{align*}
  hence the coCartesian fibration $p_0$ is classified by $\C_\o \colon  I \to \Cat$. 
  On the one hand, exponentiating and pulling back along the map \( \iota_0\colon  \Delta^0 \to I^I \) that picks the identity map we get a pullback diagram 
  \begin{align*}
    \vcenter{\xymatrix{
      **[l]  \Lax(\C_\o) \ar[d] \ar[r] \pullbackcorner & (\tint \C_\o )^I \ar[d]_{p_0^I} \ar[r] \pullbackcorner & (\tint F_\o)^I \ar[d]_{p^I} \\
      \Delta^0 \ar[r]^{\iota_0} & (I\times \Delta^{\set{0}})^I \ar[r] & (I\times \Delta^1)^I,
    }}
  \end{align*}
  On the other hand, the fiber $\E_0$ is given by the pullback diagram
  \begin{align*}
    \vcenter{\xymatrix{
      \E_0 \ar[d]_{\Lax(p)_0} \ar[r] \pullbackcorner & \E \ar[d]^{\Lax(p)} \ar[r] \pullbackcorner & (\tint F_\o)^I \ar[d]^{p^I} \\
      \Delta^{\set{0}} \ar[r] & \Delta^1 \ar[r]^-{\iota} & (I\times \Delta^1)^I.
    }}
  \end{align*}
  Observe that in both diagrams the composition of the bottom maps gives the same map \( \Delta^0 \to (I\times\Delta^1)^I \), therefore \( \Lax(\C_\o) \) and \( \E_0 \) are both pullbacks of $p^I$ along the same map, and therefore we have an equivalence \( \E_0 \simeq \Lax(\C_\o) \).
\end{proof}

\begin{mydef}
  Let \(F_\o\colon  \C_\o \to \D_\o \) be a map of $I$-diagrams of $\infty$-categories. 
  We refer to the functor 
\[ \Lax(F_\o)\colon \Lax(\C_\o) \to \Lax(\D_\o) \] of \cref{lem:assoc_lax} as the \emph{lax limit functor associated to $F_\o$}.
\end{mydef}
The associated lax limit functor has the following level-wise description.

\begin{lem} \label{lem:levelwise_assoc_lax}
  Let \(F_\o\colon  \C_\o \to \D_\o \) be a map of $I$-diagrams of $\infty$-categories and 
  \begin{align*} 
    \tilde{F} = \Lax(F_\o)\colon \Lax(\C_\o) \to \Lax(\D_\o) 
  \end{align*}
  the associated lax limit functor.
  \begin{itemize}
    \item[(1)] The associated lax limit functor $\tilde{F}$ takes a section \( x\in \Lax(\C_\o) \) to a section \( \tilde{F}(x) \in \Lax(\D_\o) \) with \( \tilde{F}(x)_a \simeq F_a(x_a) \) for every $a\in I$.
    \item[(2)] If $x$ is a coCartesian section, then so is $\tilde{F}(x)$.
  \end{itemize}
\end{lem}

\begin{proof}
  For (1) let $x$ be an object of $\Lax(\C_\o)$. 
  Consider $x$ as a vertex of $\E$ by identifying $\Lax(\C_\o)$ with $\E_0$. 
  In order to apply $\tilde{F}$ to $x$, choose a $\Lax(p)$-coCartesian edge \(e\colon  x\to y \) over \(0\to 1\). 
  Then $y$ is a vertex of $\E_1 \simeq \Lax(\D_\o)$, equivalent to $\tilde{F}(x)$.
  The edge \( e\colon  \Delta^1 \to \E \to (\tint F_\o )^I \) is $p^I$-coCartesian, and by the exponential rule defines a section \( s\colon I \times \Delta^1 \to \tint F_\o \) of $p$.
  By \cite[proposition 3.1.2.1]{Lurie} for every vertex $a$ of $I$, the induced edge \( s_a \colon  \set{a} \times \Delta^1 \to \tint F_\o \) is $p$-coCartesian.
  The restriction of the section $s$ to $I\times \Delta^{\set{0}}$ is a section of \( p_0\colon  \E_0 \to I\times\Delta^{\set{0}} \) that corresponds to \( x\colon  \Delta^0 \to \E_0 \) under the exponential rule, and the restriction of $s$ to $I\times \Delta^{\set{1}}$ corresponds to \( y\colon  \Delta^0 \to \E_1 \).
  Therefore for every $a\in I$ we have a $p$-coCartesian edge \(s_a\colon x_a \to y_a\), showing that $y_a\simeq F_a(x_a)$.
  Since $y\simeq \tilde{F}(x)$ it follows that \( \tilde{F}(x)_a \simeq y_a \simeq F_a(x_a) \) for every $a\in I$.

  For (2), assume that $x$ is a coCartesian section.
  It follows that \( s|_{I\times\Delta^{\set{0}}} \) is $p_0$-coCartesian.
  For every edge $a\to b$ in $I$ we have a commutative square
  \begin{align*}
    \vcenter{\xymatrix{
      s(0,a) \ar[r]^{s_a} \ar[d] & s(1,a) \ar[d] \\
      s(0,b) \ar[r]^{s_b} & s(1,b),
    }}
  \end{align*}
  where the horizontal arrows and the left vertical arrow are $p$-coCartesian.
  By \cite[proposition 2.4.1.7]{Lurie} the right vertical arrow \( s(a\to b,1) \) is $p$-coCartesian for every $a \to b$ in $I$, so the restriction $s|_{I\times\Delta^{\set{1}}}$ is a coCartesian section of $p_1\colon \E_1 \to I \times \Delta^{\set{1}}$, hence $\tilde{F}(x)$ is a coCartesian section.
\end{proof}

We finish this section by verifying the expected relation between $\Lax(F_\o)$ and $\invlim(F_\o)$, the induced map on the ordinary ($\infty$-categorical) limits induced by the functoriality of $\invlim$.

\begin{prop} \label{induced_functor}
  Let \(F_\o\colon  \C_\o \to \D_\o \) be a map of $I$-diagrams of $\infty$-categories. 
  The restriction of \( \Lax(F_\o)\colon  \Lax(\C_\o) \to \Lax(\D_\o) \) to $\invlim(\C_\o)$ factors through the full subcategory $\invlim(\D_\o)$ and coincides with $\invlim(F_\o)$.
\end{prop}

\begin{proof}
  The first part follows from \cref{lem:levelwise_assoc_lax}.
  Applying the exponential law to the pullback diagram \eqref{diag:Lax_as_PB} yields a commutative diagram
  \begin{align} \label{canonical_cone}
    \vcenter{\xymatrix{
      I \times \Lax(\C_\o) \ar[r] \ar[d] & \tint \C_\o \ar[d]^{q} \\
      I \ar@{=}[r] & I.
     }}
  \end{align}
  By restricting the left upper corner to $I \times \invlim(\C_\o)$, we obtain a map of coCartesian fibrations over $I$ (i.e. coCartesian edges are mapped to coCartesian edges).
  This corresponds to a map of $I$-diagrams of $\infty$-categories from the constant diagram on $\invlim(\C)$ to $\C_\o$, which is ``the universal limiting cone''.
  Let $\E \to I \times \Delta^1$ be the coCartesian fibration classified by $F_\o$ and let $\E' \subset \E$ be the full subcategory spanned by the objects of $\invlim(\C_\o)$ and $\invlim(\D_\o)$.
  Applying the same argument as before, we obtain a commutative diagram
  \begin{align} 
    \vcenter{\xymatrix{
      I \times \E' \ar[d] \ar[r] & \tint F_\o \ar[d]^{p} \\
      I \times \Delta^1 \ar@{=}[r] & I \times \Delta^1}}
  \end{align}
  which is a map of coCartesian fibrations over $I\times \Delta^1$.
  This in turn can be interpreted as a commutative square of $I$-diagrams of $\infty$-categories
  \begin{align} 
    \vcenter{\xymatrix{
      \const{\invlim(\C_\o)} \ar[d] \ar[r] & \const{\invlim(\D_\o)}\ar[d] \\
      \C_\o \ar[r] & \D_\o.
      }}
  \end{align}
It follows that the restriction of $\Lax(F_\o)$ to a functor $\invlim(\C_\o) \to \invlim(\D_\o)$ extends to a map of the universal cones and hence by the universal property of the limit must be $\invlim(F_\o)$.
\end{proof}

\begin{rem} 
  Note that the top horizontal map in diagram \eqref{canonical_cone} \emph{does not} carry coCartesian edges to coCartesian edges.
  It therefore does not correspond to a map of $I$-diagrams $\const{(\Lax(\C_\o))} \to \C_\o$.
  Instead, it should be thought of as a \emph{lax} map.
  Indeed, it is the ``universal \emph{lax} limiting cone'' on $\C_\o$.
\end{rem}

\section{Lax Limits \& Adjoint Functors} \label{sec:lax_adj}
In this section we prove an $\infty$-categorical version of a lemma of Barwick \cite{Barwick2010left} showing that for a map of $I$-diagrams of $\infty$-categories, which is level-wise left adjoint, the induced functor on the lax limits is also left adjoint.
We then use it to finish the proof of \cref{thm:main_B}.

\begin{prop}[$\infty$-categorical version of \cite{Barwick2010left}, lemma 2.23] \label{lax_adjunction}
  Let \(F_\o\colon  \C_\o \to \D_\o \) be a map of $I$-diagrams of $\infty$-categories. 
  If for all $a\in I$ the functor $F_a\colon \C_a \to \D_a$ has a right adjoint $G_a\colon \D_a \to \C_a$, 
  then the associated lax functor
  \(\Lax(F_\o)\colon  \Lax(\C_\o) \to \Lax(\D_\o) \)
  has a right adjoint \(G\colon \Lax(\D_\o) \to \Lax(\C_\o) \).
  The functor $G$ takes a section \( y_\o\in \Lax(\D_\o) \) to a section 
  \( G(y)_\o \in \Lax(\C_\o) \) with \( G(y)_a \simeq G_a(y_a) \) for every $a\in I$.
\end{prop}

\begin{proof}
  By \cref{lem:assoc_lax} the functor $\Lax(F_\o)$ is classified by the coCartesian fibration \( \Lax(p)\colon\E \to \Delta^1 \) of \cref{const:assoc_lax}.
  By \cite[definition 5.2.2.1]{Lurie} we have to show that $\Lax(p)\colon  \E \to \Delta^1$ is also a \emph{Cartesian} fibration.
  For this we need to show that for every $y\in\E$ with $\Lax(p)(y)=1$, there is a $\Lax(p)$-Cartesian edge $f\colon x\to y$ with $\Lax(p)(f)=(0\to1)$ (the unique non-degenerate edge of $\Delta^1$). 
  By \cite[proposition 2.4.1.3]{Lurie} for $f\colon x\to y$ to be $\Lax(p)$-Cartesian it suffices that its image in \( \left(\tint F_\o\right)^{I} \) is $p$-Cartesian. 
  By \cite[proposition 3.1.2.1]{Lurie} it suffices that it corresponds under the exponential rule to a map \( I\times\Delta^{1}\to\tint F_\o \) such that for all $a\in I$, 
  the edge \( \set{a} \times\Delta^{1}\to\tint F_\o\) is $p$-Cartesian. 
  We are therefore reduced to solving the lifting problem 
  \begin{align*}
    \vcenter{\xymatrix{
      I\times\Delta^{\set{1} }\ar[d]_{\Lax(p)} \ar[r] & \tint F_\o\ar[d]^{p} \\
      I\times\Delta^{1}\ar@{-->}[ru]^{h}\ar@{=}[r] & I\times\Delta^{1} ,
    }}
  \end{align*}
  such that all the edges \( h|_{\set{a} \times\Delta^{1}} \) are $p$-Cartesian in $\tint F_\o$. 
  We construct $h$ inductively over the skeleta of $I$. 
  Namely, we construct for all $n\in\N$, compatible partial lifts $h_{n}$ that make the following diagram commutative
  \begin{align*}
    \vcenter{\xymatrix{
      \sk_{n} I\times\Delta^{\set{1}} \ar[d] \ar[r] & \tint F_\o\ar[d]^{p} \\
      \sk_{n}I\times\Delta^{1}\ar@{-->}[ru]^{h_{n}}\ar@{^{(}->}[r] & I\times\Delta^{1}
    }}
  \end{align*}
  where $\sk_{n}I$ is the $n$-skeleton of $I$. 
  For $n=0$, for each $a\in I_{0}$ the pullback of $p$ along \( \set{a} \times \Delta^{1} \into I \times\Delta^{1} \) is a Cartesian fibration by the assumption that $F_a$ is a left adjoint.
  Hence, we can choose a Cartesian lift \( \set{a} \times \Delta^{1}\to\tint F_\o \) whose codomain is the image of \( \set{a} \times \Delta^{\set{1}} \).
  These Cartesian lifts assemble into a map \( h_{0}\colon\sk_{0} I \times \Delta^{1}\to\tint F_\o \) with the desired property. 
  Now assume we have constructed $h_{n}$ for some $n\ge0$. 
  For each non-degenerate simplex \( \sigma\colon\Delta^{n+1} \to I \), we can extend $h_{n}$ over $\sigma$ using the dual of \cite[proposition 2.4.1.8]{Lurie}.
  These extensions assemble into a map \( h_{n+1}\colon\sk_{n+1} I \times \Delta^{1} \to \tint F_\o \) with the desired property. 
\end{proof}
\begin{rem} \label{rem:Barwicks_lemma_fails_on_lim}
  In contrast, for \( F_\o\colon\C_{\o}\to\D_{\o} \) as in \cref{lax_adjunction}, the functor $\invlim(F_\o)$ may \emph{not} have a right adjoint. 
\end{rem}

Given a cone $\C \to \D_\o$ over an $I$-diagram of $\infty$-categories as in \cref{functorial_setting} we can expand the cone point to a constant $I$-diagram $\const{\C}$ on $\C$ and produce a map of $I$-diagrams of $\infty$-categories $F_\o \colon \const{\C} \to \D_\o$, where $F_a\colon\C \to \D_a$ are the functors of the given cone. 
More formally, the cone $\C \to \D_\o$ is given by a map $I^\cone \to \Cat$. 
There is a map of simplicial sets $I \times \Delta^1 \to \Delta^0 \star I$ obtained as the composition of the canonical maps (see \cite[section 4.2.1]{Lurie})
\[ I \times \Delta^1 \to \Delta^0 \diamond I \to \Delta^0 \star I = I^\cone \]
and $F_\o \colon \const{\C} \to \D_\o$ corresponds to the composition $I \times \Delta^1 \to I^\cone \to \Cat$. 
Combining this with \cref{Const_diag}, we obtain an associated lax limit functor $\tilde{F} = \Lax(F_\o)\colon \C^I \to \Lax(\D_\o)$, which takes an $I$-diagram $x_\o \in \C^I$ to a section $\tilde{F}(x)_\o \in \Lax(\D_\o)$ with $\tilde{F}(x)_a \simeq F_a(x_a)$.
Applying \cref{lax_adjunction} we obtain:
\begin{cor} \label{cor:const_lax_adjunction}
  In the \cref{functorial_setting}, if for every $a\in I$ the functor \( F_a\colon\C \to \D_a \) of the cone $\C\to \D_\o$ has a right adjoint \( G_a\colon\D_a\to \C_a\), then the associated functor \( \tilde{F}\colon\C^I \to \Lax(\D_\o) \) has a right adjoint 
  \begin{align*}
    \tilde{F}\colon\C^I \adj \Lax(\D_\o)\noloc\tilde{G}^I ,
  \end{align*} 
  which takes a section \(y_\o\in\Lax(\D_\o) \) to an $I$-diagram $\tilde{G}^I(y)_\o\in \C^I$ with \(\tilde{G}^I(y)_a \simeq G_a(y_a) \).
\end{cor}

Before we proceed, we need the following easy lemma:
\begin{lem}[Restriction of an adjunction] \label{lem:restricted_adjunction}
  Let \( F\colon\C \adj \D\noloc G \) be an adjunction of $\infty$-categories and $\D' \subset \D$ a full subcategory. 
  If $F$ factors via $\D'$ then \( F\colon\C \to D' \) is left adjoint, with right adjoint given by the restriction \( G|_{\D'} \colon \D' \to \C \).
\end{lem}
\begin{proof}
  By \cite[proposition 5.2.2.8]{Lurie} there exists a unit transformation \(u\colon Id_\C \to G \circ F \).
  Because \( G \circ F = G|_{\D'} \circ F \) we can consider $u$ as a morphism \( u \colon Id_\C \to G|_{\D'} \circ F \).
  To see that $u$ is a unit transformation for \( F\colon\C \adj \D' \noloc G|_{\D'} \), note that for every $x\in \C, \, y\in \D'$ we have \( \Map_{\D'}(Fx,y) \simeq \Map_\D(Fx,y) \), since $\D'\subset \D$ is a full subcategory.
\end{proof}

We are now ready to prove our main result, which is a more detailed version of \cref{thm:main_B}.
\begin{thm} \label{inf_adjoint_descent}
  Let $I$ be a simplicial set, $\D_\o$ an $I$-diagram of $\infty$-categories, and $\C \to \D_\o$ a cone.
  Assume that for every $a\in I$ the functor \( F_a\colon \C \to \D_a \) of the cone $\C\to \D_\o$ has a right adjoint \( G_a\colon \D_a\to \C_a\) and $\C$ has all $I$-limits.
  The comparison functor 
  \( F\colon \C \to \invlim \D_\o \) 
  has a right adjoint $G$, given as a composition
  \( \invlim \D_\o \xto{G^I} \C^I \xto{\invlim} \C \), 
  where $G^I$ is the restriction of the functor \( \tilde{G}^I\colon \Lax(\D) \to \C^I \) of \cref{cor:const_lax_adjunction} to the full subcategory \( \invlim(\D_\o) \subset \Lax(\D_\o) \).
\end{thm}
\begin{proof}
  Let \( \tilde{F} \colon \C^I \adj \Lax(\D_\o)\noloc \tilde{G}^I \) be the adjunction of \cref{cor:const_lax_adjunction}. 
  The diagonal functor  \(\Delta\colon  \C \to \C^I\) also has a right adjoint \( \invlim \colon  \C^I \to \C \).
  By \cite[proposition 5.2.2.6]{Lurie} the composition \( \C \xto{\Delta} \C^I \xto{\tilde{F}} \Lax(\D) \) has a right adjoint given by the composition \( \Lax(\D) \xto{\tilde{G}^I} \C^I \xto{\invlim} \C \). 

  By \cref{induced_functor} we get a commutative square
  \begin{align*}
    \vcenter{\xymatrix{
      \invlim (\const{\C}) \ar[d]_{\invlim(F_\o)} \ar@{^(->}[r] & \C^I \ar[d]^{\tilde{F}} \\
      \invlim (\D_\o) \ar@{^(->}[r] & \Lax(\D_\o) . \\
    }}
  \end{align*}
  Since the diagonal functor \( \Delta\colon \C \to \C^I \) factors through \( \invlim (\const{\C}) \), the composition $\tilde{F}\circ \Delta$ factors through the full subcategory \( \invlim (\D_\o) \subset \Lax(\D_\o) \).
  \Cref{lem:restricted_adjunction} shows that \( \C \to \invlim \D_\o \) has a right adjoint, given by the composition
  \( \invlim \D_\o \xto{G^I} \C^I \xto{\invlim} \C . \)
  Finally, we observe that the left adjoint \( \C \to \invlim \D_\o \) being the composition $\C \to \invlim (\const{\C}) \to \invlim(\D_\o)$ is just the comparison functor of the cone \( \C \to \D_\o \).
\end{proof}

\bibliography{references}
\bibliographystyle{alpha}

\end{document}

%% file: preamble.tex
\usepackage{amssymb}
\usepackage{amsfonts}
\usepackage{graphicx}
\usepackage{amstext}

\usepackage{amsmath,amscd, amsthm}

\usepackage{euscript}
\usepackage[latin9]{inputenc}
\usepackage{geometry}
\usepackage{pdfsync}

\usepackage{hyperref}

\usepackage[all]{xy}

\usepackage{stmaryrd}

\usepackage{tikz} 
\usetikzlibrary{arrows,decorations.pathmorphing,backgrounds,positioning,fit,petri,shapes.misc,decorations.markings,shapes,shapes.multipart,calc}

\usepackage{cleveref}

\tikzset{
	plain trans/.style = {circle, minimum size = 14, inner sep=2, draw, font=\scriptsize}}
\tikzset{
	plain inv trans/.style = {circle, minimum size = 18, inner sep=1, draw, font=\tiny}}
\tikzset{
	plain dual trans/.style = {circle, minimum size = 14, inner sep=2, draw, font=\scriptsize}}
\tikzset{
	plain inv dual trans/.style = {circle, minimum size = 18, inner sep=1, draw, font=\tiny}}

\tikzset{trans/.style={circle, minimum size = 14, inner sep=1, draw, font=\scriptsize,
	path picture={%
	\filldraw[anchor=center,black]
	($(path picture bounding box.south east)!0.15!(path picture bounding box.north east)$) to [out=-135, in=0]	(path picture bounding box.south) to [out=180, in=-45]
	($(path picture bounding box.south west)!0.15!(path picture bounding box.north west)$) to cycle;
      }}}
      
\tikzset{inv trans/.style={circle, minimum size = 14, inner sep=1, draw, font=\tiny,
	path picture={%
	\filldraw[anchor=center,black]
	($(path picture bounding box.south east)!0.15!(path picture bounding box.north east)$) to [out=-135, in=0]	(path picture bounding box.south) to [out=180, in=-45]
	($(path picture bounding box.south west)!0.15!(path picture bounding box.north west)$) to cycle;
      }}}      

\tikzset{dual trans/.style={circle, minimum size = 14, inner sep=1, draw, font=\scriptsize,
	path picture={%
	\filldraw[anchor=center,black]
	($(path picture bounding box.north east)!0.15!(path picture bounding box.south east)$) to [out=135, in=0]	(path picture bounding box.north) to [out=180, in=45]
	($(path picture bounding box.north west)!0.15!(path picture bounding box.south west)$) to cycle;
      }}}      

\tikzset{inv dual trans/.style={circle, minimum size = 14, inner sep=1, draw, font=\tiny,
	path picture={%
	\filldraw[anchor=center,black]
	($(path picture bounding box.north east)!0.15!(path picture bounding box.south east)$) to [out=135, in=0]	(path picture bounding box.north) to [out=180, in=45]
	($(path picture bounding box.north west)!0.15!(path picture bounding box.south west)$) to cycle;
      }}}      

\tikzset{
	func/.style = {font=\scriptsize},
	in func/.style = {font=\scriptsize, above, text opacity=0},
	in func visible/.style = {font=\scriptsize, above},
	out func/.style = {font=\scriptsize, below, text opacity=0},	
	out func visible/.style = {font=\scriptsize, below},	
	mid func/.style = {font=\scriptsize, inner sep=1, fill=white} 
}

\tikzset{
	cat/.style = {font=\scriptsize}}
\tikzset{
	frame/.style = {double,rounded corners}}

\tikzset{
	paste/.style = {dashed,opacity=0.5}}

\tikzset{
    dir/.style={},
    dir end/.style={},    
    dual/.style={},
    eq/.style={inner sep=0pt},
    plain/.style={} 
}

\newcommand\noloc{%
  \nobreak
  \mspace{6mu plus 1mu}
  {:}
  \nonscript\mkern-\thinmuskip
  \mathpunct{}
  \mspace{3mu}
}
\renewcommand{\o}{\bullet}
\newcommand{\N}{\mathbb{N}}
\newcommand{\Z}{\mathbb{Z}}

\newcommand{\Map}{\operatorname{Map}}

\newcommand{\sk}{\operatorname{sk}}
\newcommand{\iso}{\overset{\sim}{\longrightarrow}}
\renewcommand{\=}{\cong}

\newcommand{\Aut}{\operatorname{Aut}}
\newcommand{\BAut}{\operatorname{BAut}}
\newcommand{\hAut}{\operatorname{hAut}}
\newcommand{\BMap}{\operatorname{BMap}}
\newcommand{\Conj}{\operatorname{Conj}}

\def\invlim{\qopname\relax m{\underleftarrow{\lim}}}
\def\colim{\qopname\relax m{\underrightarrow{\lim}}}

\newcommand{\Lax}{\underleftarrow{\operatorname{Lax}}}
\newcommand{\const}[1]{#1^{\text{ct}}}

\newcommand{\cone}{\triangleleft}

\newcommand{\into}{\hookrightarrow}

\newcommand{\xto}{\xrightarrow}

\newcommand{\set}[1]{\left\{{#1}\right\}}

\newcommand{\tint}{\textstyle\int}

\newcommand{\E}{\mathcal{E}}

\newcommand{\C}{\mathcal{C}}
\newcommand{\D}{\mathcal{D}}

\newcommand{\U}{\mathcal{U}}

\newcommand{\vect}{\mathbf{Vect}}

\newcommand{\Cat}{\mathbf{Cat}_{\infty}}
\newcommand{\adj}{\leftrightarrows}

\newtheorem{thm}{Theorem}[section]
\newtheorem{prop}[thm]{Proposition}
\newtheorem{cor}[thm]{Corollary}
\newtheorem{lem}[thm]{Lemma}

\theoremstyle{definition}
\newtheorem{mydef}[thm]{Definition}
\newtheorem{rem}[thm]{Remark}
\newtheorem{ex}[thm]{Example}
\newtheorem{notation}[thm]{Notation}
\newtheorem{construction}[thm]{Construction}
\newtheorem{mysetting}[thm]{Setting}
\newtheorem{theorem}{Theorem}

\newcommand{\pullbackcorner}[1][dr]{\save*!/#1-1.2pc/#1:(1,-1)@^{|-}\restore}

